\RequirePackage{rotating}
\documentclass[a4paper,twopage,reqno,12pt]{amsart}
\usepackage{eurosym}
\usepackage[top=30mm,right=30mm,bottom=30mm,left=30mm]{geometry}
\usepackage{tabularx}
\usepackage{graphicx}
\usepackage{longtable}

\usepackage{amsmath}
\usepackage{amssymb}
\usepackage{amsfonts}
\usepackage{amsthm}
\usepackage{amstext}
\usepackage{amsbsy}
\usepackage{amsopn}
\usepackage{amscd}
\usepackage{enumerate}
\usepackage{booktabs}
\usepackage[caption=false]{subfig}
\usepackage{color}
\usepackage[colorlinks]{hyperref}
\usepackage[hyperpageref]{backref}
\usepackage{multirow}
\usepackage{longtable}
\usepackage{float}

\usepackage{lscape}
\usepackage{pdflscape}
\usepackage{url}
\usepackage{adjustbox}
\usepackage{rotating}
\usepackage{etoolbox}
\usepackage{xcolor}
\usepackage{comment}
\newtheorem{theorem}{Theorem}[section]

\newtheorem{lemma}[theorem]{Lemma}
\newtheorem{proposition}[theorem]{Proposition}

\AtEndEnvironment{proof}{\setcounter{claim}{0}}
\theoremstyle{definition}

\newtheorem{example}[theorem]{Example}

\numberwithin{equation}{section}

\newcommand{\PSL}{\mathrm{PSL}}

\newcommand{\PGL}{\mathrm{PGL}}

\newcommand{\Out}{\mathrm{Out}}

\newcommand{\PG}{\mathrm{PG}}
\newcommand{\AG}{\mathrm{AG}}

\newcommand{\Dmc}{\mathcal{D}}

\renewcommand{\leq}{\leqslant}
\renewcommand{\geq}{\geqslant}

\begin{document}
\title[Flag-transitive $2$-designs with $\lambda =3$] {$2$-$(v,k,3)$ designs admitting an almost simple, flag-transitive automorphism group with socle $\PSL(2,q)$}
\author[Hongxue Liang]{Hongxue Liang}
\address{School of Mathematics, Foshan University, Foshan 528000, P.R. China}
\email{hongxueliang@fosu.edu.cn}
\author[Zhihui Liu]{Zhihui Liu}
\address{School of Mathematics, Foshan University, Foshan 528000, P.R. China}
\email{zhihuiliu889@gmail.com}
\author[Alessandro Montinaro]{Alessandro Montinaro}
\thanks{Corresponding author: Alessandro Montinaro}
\address{Dipartimento di Matematica e Fisica “E. De Giorgi”, University of Salento, Lecce, Italy}
\email{alessandro.montinaro@unisalento.it}
\subjclass[MSC 2020:]{05B05; 05B25; 20B25}
\date{\today }

\begin{abstract}

In this paper, we completely classify the non-trivial $2$-$(v,k,3)$ designs admitting an almost simple, flag-transitive automorphism group with socle $\PSL(2,q)$. 

\end{abstract}

\maketitle

\section{Introduction and Main Result } \label{Intro}
A $2$-$(v,k,\lambda )$ design $\Dmc$ is a pair $(\mathcal{P},%
\mathcal{B})$ with a set $\mathcal{P}$ of $v$ points and a set $\mathcal{B}$
of $b$ blocks such that each block is a $k$-subset of $\mathcal{P}$ and each two distinct points are contained in $\lambda $ blocks. The \emph{replication number} $r$ of $\mathcal{D}$ is the number of blocks containing a given point.
We say $\mathcal{D}$ is \emph{non-trivial} if $2<k<v-1$, and \emph{symmetric} if $v=b$. All $2$-$(v,k,\lambda
)$ designs in this paper are assumed to be non-trivial.

An automorphism of $%
\mathcal{D}$ is a permutation of the point set which preserves the block
set. The set of all automorphisms of $\mathcal{D}$ with the composition of
permutations forms a group, denoted by $\mathrm{Aut(\mathcal{D})}$. For a
subgroup $G$ of $\mathrm{Aut(\mathcal{D})}$ acting point-transitively on $\mathcal{D}$, $G$ is said to be \emph{%
point-primitive} if $G$ acts primitively on $\mathcal{P}$, and said to be 
\emph{point-imprimitive} otherwise. In this setting, we also say that $%
\mathcal{D}$ is either \emph{point-primitive} or \emph{point-imprimitive}, respectively. A \emph{flag} of $\mathcal{D}$ is a pair $(x,B)$ where $x$ is
a point and $B$ is a block containing $x$. If $G\leq \mathrm{Aut(\mathcal{D})%
}$ acts transitively on the set of flags of $\mathcal{D}$, then we say that $%
G$ is \emph{flag-transitive} and that $\mathcal{D}$ is a \emph{flag-transitive design}.

There is a vast literature focusing on $2$-designs $\mathcal{D}$ admitting a
flag-transitive automorphism group $G$. If $\lambda =1$, then $G$ acts
point-primitively on $\mathcal{D}$ by a\ famous result by Higman-McLaughlin 
\cite[Proposition 3]{HM}. However, the conclusion of the Higman-McLaughlin
result is no longer true for $\lambda >1$, and there are several well-known
examples of flag-transitive point-imprimitive $2$-designs. Nevertheless, $G$
is still point-primitive for $\lambda \le 4$ except for eleven numerical
cases by a recent result by Devillers-Praeger \cite[Theorem 1.1]{DevillersPraeger23}. The
importance of the point-primitivity relies on the fact that the Theorem of
O'Nan-Scott can be applied to determine the general structure of $G$.
Indeed, using the previous result, and as a summary of various papers due to
several authors, a satisfactory classification of $(\mathcal{D},G)$ was
achieved by Buekenhout et al. in \cite{ Buekenhout90} for $\lambda =1$, and by the
first and the third author in \cite{LM} for $\lambda =2$. Hence, $\lambda =3$
is the next natural step to analyze, and in this case, by \cite[Theorem 1.1]
{DevillersPraeger23},  the unique flag-transitive, point-imprimitive $2$-design with $
\lambda =3$ is the $2$-$(45,12,3)$ design constructed by Praeger \cite{P}. 

Our paper is a contribution to the classification of the flag-transitive $2$
-designs with $\lambda =3$. In particular, we completely classify the pair $(
\mathcal{D},G)$ when $\lambda =3$ and $\PSL(2,q)\trianglelefteq G\leq \mathrm{P\Gamma
L(2,q)}$, the collineation group of the projective line $\PG(1,q)$. The result
contained in our paper is meant as a natural prosecution of the analysis
involving $\PSL(2,q)\trianglelefteq G\leq \mathrm{P\Gamma L(2,q)}$ and carried out by
Delandtsheer \cite{De} for $\lambda =1$, by Montinaro, Zhao, Zhang and Zhou 
\cite{MZZZ} for $\lambda =2$, and by Alavi, Bayat and Daneshkhah \cite{AlaviSymPSL2}
for any value $\lambda $ when $\mathcal{D}$ is symmetric. More precisely, we prove the following result.

\begin{theorem}\label{main}
Let $\mathcal{D}$ be a non-trivial $2$-design with $\lambda=3$ 
admitting $\PSL(2,q) \unlhd G \leq \mathrm{P \Gamma L(2,q)}$ as a flag-transitive automorphism group. Then one of the following holds:
\begin{enumerate}
    \item $\mathcal{D}$ is the complete $2$-$(5,3,3)$ design and $\PSL(2,4)\unlhd G \leq \mathrm{P \Gamma L(2,4)}$;
    \item $\mathcal{D}$ is the $2$-$(8,4,3)$ design isomorphic to $\AG(3,2)$ and $G \cong \PSL(2,7)$;
    \item $\mathcal{D}$ is the $2$-$(11,3,3)$ design as in Example \ref{Ex1}, and $G \cong \PSL(2,11)$;
    \item $\mathcal{D}$ is the $2$-$(11,6,3)$ design, complement of the $2$-$(11,5,2)$ Paley design, and $G \cong \PSL(2,11)$;
    \item $\mathcal{D}$ is the $2$-$(26,6,3)$ design as in Example \ref{Ex2} and $\PSL(2,25) \unlhd G \leq \mathrm{P \Sigma L(2,25)}$.  
\end{enumerate}
\end{theorem}

\bigskip

The result is obtained as follows: we use \cite[Theorem 1.1]{DevillersPraeger23} to prove that $G$ acts point-primitively on $\mathcal{D}$, and that no novelties occur by \cite[Theorems 1.1 and 2.2]{gmaxi}, that is, also $Soc(G)=\PSL(2,q)$ acts point-primitively on $\mathcal{D}$. Finally, we complete the proof by combining the list of the subgroups of $\PSL(2,q)$ provided in \cite{Di, Hup} with some constraints on the parameters $\mathcal{D}$ arising from the condition $\lambda=3$ and some remarkable geometric configurations left invariant by $G$. 

\bigskip

The following examples correspond to cases (3) and (5), respectively. In both cases $G$ acts point-$2$-transitively on $\mathcal{D}$ and (5) belongs to an infinite family of examples with $q$ a square and $\lambda=q^{1/2}+1$ when $q=25$.

\bigskip

\begin{example}\label{Ex1}
Let $G=PSL(2,11)$ and let $P,Q$ be subgroups of $G$ isomorphic to $A_{5}$ and $D_{12}$, respectively, such that $P \cap Q=S_{3}$, and let $\mathcal{D}=\left(\mathcal{P},\mathcal{B}\right)$, where $\mathcal{P}=\{Px:x\in G\}$ and $\mathcal{B}=\{Qx:x\in G\}$. Then $\mathcal{D}$ is a $2$-$(11,3,3)$ design admitting $G$ as a full flag-transitive, point-$2$-transitive automorphism group. 

\end{example}
\begin{proof}
The incidence structure $\mathcal{D}$ has parameters $(v,k,b,r)=(11,3,55,15)$ and has $G$ as a flag-transitive automorphism group by \cite[Lemmas 1 and 2]{HM}. Actually, $\mathcal{D}$ is a $2$-$(11,3,\lambda)$ design since $G$ acts point-$2$-transitively on $\mathcal{D}$, and $\lambda=\frac{r(k-1)}{v-1}=3$.

Let $A=Aut(\mathcal{D})$. Clearly, one has $G \leq A$. Suppose that $A \neq G$. Then either $A=S_{11}$ or $M_{11}$ by \cite[Table B.4]{DM}. Then $b=\binom{11}{3}=165$ since in both cases $A$ acts point-$3$-transitively on $\mathcal{D}$, a contradiction. Thus $A=G$. This completes the proof.
\end{proof}

\begin{example}\label{Ex2}
The incidence structure having the projective line $\PG(1,25)$ as a set of points, and one of the two orbits under $\PSL(2,5)\unlhd G \leq \mathrm{P \Sigma L(2,25)}$ of the Baer sublines $\PG(1,5)$ of $\PG(1,25)$ is a $2$-$(26,6,3)$ design admitting $G$ as a flag-transitive point-transitive  automorphism group.
\end{example}

\begin{proof}
Consider the natural $2$-transitive action of $\PGL(2,q)$ on the projective line $\PG(1,q)$. If $q$ is a square, the Baer sublines $\PG(1,q^{1/2})$ are projectively equivalent, that is to say, $\PGL(2,q)$ acts transitively on the set of the Baer sublines $\PG(1,q^{1/2})$ inducing a copy of $\PGL(2,q^{1/2})$ on each of them. Now, the incidence structure having $\mathcal{P}=\PG(1,q)$ as a set of points, and the set $\mathcal{B}$ of whole Baer sublines of $\PG(1,q)$ as a block set, is a $2$-$(q+1,q^{1/2}+1,q^{1/2}+1)$ design since $G$ acts $2$-transitively on $\PG(1,q)$ and $\left\vert \PGL(2,q):\PSL(2,q^{1/2})\right\vert=q^{2}(q+1)$. Now, $\PSL(2,q)$ acts $2$-transitively on $\PG(1,q)$ and partitions the set $\mathcal{B}$ of Baer sublines of $\PG(1,q)$ into $\gcd(2,q-1)$ sets, say $\mathcal{B}_{1},...,\mathcal{B}_{\gcd(2,q-1)}$, both of size $q^{2}(q+1)/2$.
Hence, the incidence structures $\mathcal{D}_{i}=(\mathcal{P},\mathcal{B}_{i})$, $i=1,...,\gcd(2,q-1)$, are $\gcd(2,q-1)$ isomorphic $2$-$\left(q+1,q^{1/2}+1,\frac{q^{1/2}+1}{\gcd(2,q-1)}\right )$ designs under an element $\mathrm{P \Gamma L(2,q)}\setminus\mathrm{ P\Sigma L(2,q)}$. Moreover, since $\PSL(2,q)$ induces a copy of $\PGL(2,q^{1/2})$ on each of $\mathcal{B}_{i}$, the group $\PSL(2,q) \unlhd G \leq \mathrm{P \Sigma L}(2,q)$ acts flag-transitively on the $2$-designs $\mathcal{D}_{i}$, which have $\lambda=3$ for $q=25$. 
\end{proof}

\section{Preliminaries}

In this section, we provide some useful facts in both design theory and group theory.
If a group $G$ acts transitively on a set $\mathcal{P}$ and $\alpha\in \mathcal{P}$, then the
\emph{subdegree} $d$ of $G$ is the length of some $G_\alpha$-orbit.  We say that $d$ is non-trivial if the orbit is not $\{\alpha\}$.

\begin{lemma} \label{lem:basic-params}
Let $\mathcal{D}$ be a $2$-$(v,k,3)$ design that admits a flag-transitive
automorphism group $G$. Let $\alpha\in \mathcal{P}$ and $G_\alpha$ be the point stabilizer, then
\begin{enumerate}
\item[\rm(i)]   $r(k-1)=3(v-1)$ and $vr=bk$;
\item[\rm(ii)]  $b\ge v$ and $r \ge k$;
\item[\rm(iii)]  $r \mid |G_\alpha|$ and $r^{2}>3v$;
\item[\rm(iv)]  $r\mid 3d$ for every nontrivial subdegree $d$ of $G$.
\end{enumerate}
\end{lemma}

\begin{lemma} \label{lem:dp}
Let $\mathcal{D}$ be a $2$-$(v,k,3)$ design that admits a flag-transitive
automorphism group $G$. Then one of the following holds:
\begin{enumerate}
\item[\rm(i)]   $G$ acts point-primitively on $\mathcal{D}$;
\item[\rm(ii)]  $\mathcal{D}$ is the $2$-$(45,12,3)$ design constructed in \cite{P}, and $G$ is one of the groups $(3^{4}:5):8$, $3^{4}:2.A_{5}$ or $3^{4}:2.S_{5}$ acting point-imprimitively on $\mathcal{D}$.
\end{enumerate}
\end{lemma}

\begin{proof}
The assertion immediately follows from \cite[Theorem 1.1]{DevillersPraeger23}.
\end{proof}

\bigskip
On the basis of Lemma \ref{lem:dp}, in the sequel we may assume that $G$ acts point-primitively on $\mathcal{D}$. 
\bigskip




\begin{lemma}\label{lem:HX}
{\rm (\cite[Lemma~2.2]{AlaviSymPSL2})}
Let $G$ be an almost simple group with socle $X$, and let $H$ be maximal
in $G$ not containing $X$. Then $G=HX$, and $|H|$ divides $|\Out(X)|\cdot |H\cap X|$.
\end{lemma}

\begin{lemma}\label{subdegrees-dihedral}
{\rm (\cite[Table~2]{FaradzevIvanov})}
Let $G=\PSL(2,q)$ act on the set of cosets of its
subgroup $H=\mathrm{D}_{2(q-\epsilon)/\gcd(2,q-1)}$ 
with $\epsilon=\pm1$, respectively. Then the subdegrees of $G$ are as given in
Table~\ref{tab:subdegrees-PSL2-dihedral}, where $a^{b}$ means that the subdegree $a$
appears with multiplicity $b$.
\end{lemma}

\begin{table}[htbp]
\centering

\caption{The subdegrees of $\PSL(2,q)$ on the set of cosets of its dihedral subgroup $H$.}
\label{tab:subdegrees-PSL2-dihedral}
\begin{tabular}{llll}
\hline
$G$ & $H$ & subdegrees & Conditions \\
\hline
$\PSL(2,q)$ & $\mathrm{D}_{2(q-1)}$   & $1,\ (q-1)^{q/2-1},\ 2(q-1)$  & $q$ even \\
$\PSL(2,q)$ & $\mathrm{D}_{2(q+1)}$   & $1,\ (q+1)^{q/2-1}$  & $q$ even \\
$\PSL(2,q)$ & $\mathrm{D}_{(q\pm1)}$  & $1,\ \bigl(\tfrac{q\pm1}{4}\bigr)^{2},\
     \bigl(\tfrac{q\pm1}{2}\bigr)^{(q\pm1)/2-2},\
     (q\pm1)^{(q+2\mp5)/4}$
  & $q\equiv\pm3 \pmod 8$ \\
$\PSL(2,q)$ & $\mathrm{D}_{(q\pm1)}$
  & $1,\ \bigl(\tfrac{q\pm1}{4}\bigr)^{2},\
     \bigl(\tfrac{q\pm1}{2}\bigr)^{(q\pm1)/2-2},\
     (q\pm1)^{(q+2\mp5)/4}$
  & $q\equiv\mp1 \pmod 8$ \\
\hline
\end{tabular}
\end{table}


\section{Proof of the main result}

In this section, we prove Theorem~1.1.
Let $\mathcal{D}$  be a $2$-$(v,k,3)$ design that admits a flag-transitive an almost simple automorphism group $G$ with socle $X=\mathrm{PSL}(2,q)$. Then $G$ acts point-primitively on $\mathcal{D}$ by Lemma \ref{lem:dp}. Note that $\mathcal{D}$ is the $2$-$(8,4,3)$ design, isomorphic to $\AG(3,2)$, and $G=\PSL(2,7)$ when $\gcd(r,\lambda)=1$ by \cite{rlambdacoprime}, and the $2$-$(11,6,3)$ design, complementary design of the $2$-$(11,5,2)$ Paley design, and $G=\PSL(2,11)$ when $\mathcal{D}$ is symmetric \cite{AlaviSymPSL2}. These are cases (2) and (4) of Theorem \ref{main} respectively. Therefore, in the following, we focus on the case where $\mathcal{D}$ is a non-symmetric $2$-$(v,k,3)$ design with $3\mid r$.

\bigskip

\begin{lemma}\label{lem:Xprim}
$X$ acts point-primitively on $\mathcal{D}$.    
\end{lemma}

\begin{proof}
Since $G$ is point-primitive on $\mathcal{D}$ and $X \unlhd G$, it follows that $X$ acts point-transitively on $\mathcal{D}$. Then, by \cite{gmaxi}, either $X_\alpha$ is the maximal subgroup of $X$, and hence $X$ acts point-primitively on $\mathcal{D}$ , or $(G,G_\alpha,|G:G_\alpha|,R)$ is as in Table \ref{Tab1}. 

\begin{table}[h!]
\footnotesize
\caption{The possibilities of  $(G,G_{\alpha},|G:G_\alpha|,R)$.}
\label{Tab1}
\footnotesize
\begin{tabular}{lllll}
\hline
Line & $G$ & $G_\alpha$ & $|G:G_\alpha|$ & $R$ \\
\hline
1 & $\mathrm{PGL}(2,7)$         & $\mathrm{D_{12}}$                & 28 &  3\\  
2 & $\mathrm{PGL}(2,7)$         & $\mathrm{D_{16}}$                & 21 &  4\\  
3 & $\mathrm{PGL}(2,9)$         & $\mathrm{D_{20}}$                & 36 &  5\\  
4 & $\mathrm{PGL}(2,9)$         & $\mathrm{D_{16}}$                & 45 &  4\\  
5 & $\mathrm{M_{10}}$           & $C_{5}\rtimes C_4$               & 36 &  5\\  
6 & $\mathrm{M_{10}}$           & $C_{8}\rtimes C_{2}$             & 45 &  4\\  
7 & $\mathrm{P\Gamma L}(2,9)$   & $C_{10}\rtimes C_{4}$            & 36 &  5\\  
8 & $\mathrm{P\Gamma L}(2,9)$   & $C_{8}\cdot Aut(C_8)$            & 45 &  4\\ 
9 & $\mathrm{PGL}(2,11)$        & $\mathrm{D_{20}}$                & 66 &  5\\ 
10& $\mathrm{PGL}(2,q),q=p\equiv \pm 11,19(\text{mod}~ 40)$ & $\mathrm{S_{4}}$  &$\frac{q(q^2-1)}{24}$  &  6\\
\hline
\end{tabular}
\end{table}
Suppose that $X_\alpha$ is not maximal in $X$. Then one of the lines of Table \ref{Tab1}
holds. As $G$ is point-primitive, we have $v=|G:G_\alpha|$ as in the fourth column of 
Table \ref{Tab1}. Let $R:=(|G_\alpha|,3(v-1))$. We list $R$ in the fifth column of Table \ref{Tab1}. As $r\mid R$ and $r^2 >3v$, we have $R^2>3v$, ruling out all the cases in Lines 1-9. For the case as in Line 10, $(G,G_\alpha,v)=(\PGL(2,q), S_4,\frac{q(q^2-1)}{24})$, by $r\leq |G_{\alpha}|$ and $r^2>3v$, we have $p\leq 16$. Combining with $p\equiv \pm 11,19(\text{mod}~ 40)$, we obtain $p=11$, and so $v=55$, $R=6$, contradicting the fact that $R^2>3v$. This completes the proof.
\end{proof}

\bigskip

It follows from Lemma \ref{lem:Xprim} that $X_\alpha$ is the maximal subgroup of $X$. Hence, by \cite{Di, Hup}, $X_\alpha$ is one of the following groups:

\begin{enumerate}
  \item $\mathrm{PGL}(2,q_0)$, for $q=q_0^2$ odd;
  \item $\mathrm{PSL}(2,q_0)$, for $q=q_0^t$ odd, where $t$ is odd and prime;
  \item $\mathrm{PGL}(2,q_0)$, for $q=2^f=q_0^t$, where $t$ is prime and $q_0\neq 2$;
  \item $\mathrm{A_4}$, for $q=p\equiv \pm 3 \text{(mod 8)}$ and $q\not\equiv \pm 1\text{(mod 10)}$;
  \item $\mathrm{S_4}$, for $q=p\equiv \pm 1\text{(mod 8)}$;
  \item $\mathrm{A_5}$, for $q\equiv \pm 1\text{(mod 10)}$, where either $q=p$ or $q=p^2$ and $p\equiv \pm 3\text{(mod 10)}$;
  \item $\mathrm{D}_{2(q-1)/\gcd(2,q-1)}$;
  \item $\mathrm{D}_{2(q+1)/\gcd(2,q-1)}$;
  \item $C_p^f\rtimes C_{q-1/\gcd(2,q-1)}$.  
\end{enumerate}   
In what follows, we will analyze each of these cases separately.

\begin{lemma}\label{case1}
$X_{\alpha}\neq \mathrm{PGL}(2,q_0)$, for $q=q_0^2$ odd.
\end{lemma}

\begin{proof}
Suppose that $X_{\alpha}=\mathrm{PGL}(2,q_0)$, for $q=q_0^2$ odd.
Then we obtain $v=\frac{q_0(q_0^2+1)}{2}$.
It follows from Lemma~\ref{lem:basic-params}\,(iv) and \cite[Theorem~2]{CameronMaimaniOmidi} that
$r\mid 3q_0(q_0^2-1)$. 
Combining this with $r\mid 3(v-1)$, $\gcd (v-1,q_0)=1$ and $\gcd(q_0+1,q_0^2+q_0+2)=2$, 
we conclude that $r \mid 6(q_0-1)$.
From $r^2>3v$, we obtain $q_0(q_0^2+1) < 24(q_0-1)^2$, implying 
$q_0=5,7,9,11,13,17,19$, and so $v=65,175,369,671,1105,2465,3439$, respectively.
Since $r\mid 6(q_0-1)$ and $r^2>3v$, the corresponding values of $r$ are
$24,36,48,60,72,96,108$, respectively. Recall that
$r(k-1)=3(v-1)$, $b=vr/k$ is a positive integer, and $r>k$.
Thus, $(v,b,r,k,\lambda)=(369,738,48,24,3)$, $X \cong \PSL(2,3^{4})$, $X_{\alpha}=\PGL(2,3^{2})$. Moreover, $\PSL(2,9) \leq X_{B}$ since the order of $X_{B}$ is divisible by $45$. Thus either $X_{B}=\PSL(2,9)$ and $X$ acts block-transitively on $\mathcal{D}$, or $X_{B}=\PGL(2,9)$ and $X$ partitions the block-set of $\mathcal{D}$ into two orbits of length $369$. 

Let $H$ be a Sylow $5$-subgroup of $X$. Then $N_{X}(H)\cong \mathrm{D_{80}}$ and $N_{X_{\alpha}}(H)\cong \mathrm{D_{20}}$, and hence $H$ fixes exactly $4$ points on $\mathcal{D}$. On the other hand, either $X$ acts block-transitively on $\mathcal{D}$ and $N_{X_{B}}(H)\cong \mathrm{D_{10}}$, or $X$ partitions the block-set of $\mathcal{D}$ into two orbits of length $369$ and $N_{X_{B}}(H)\cong \mathrm{D_{20}}$. However, in either case $H$ leaves invariant exactly $8$ blocks of $\mathcal{D}$. Now, if $X$ acts block-transitively on $\mathcal{D}$ either $B$ consists of two $X_{B}$-orbits each of length $12$, or four $X_{B}$-orbits each of length $6$; if $X$ partitions the block-set of $\mathcal{D}$ into two orbits of length $369$, the group $X_{B}=\PSL(2,9)$ partitions $B$ into four orbits of length $6$. Thus, the number of fixed points by $H$ on each $H$-invariant block is constant and is either $2$ or $4$. In the latter case, any of the blocks left invariant by $H$ contains the $4$ points fixed by $H$, forcing the number of blocks fixed by $H$ to be $3$ since $\lambda=3$, whereas $H$ leaves invariant exactly $8$ blocks of $\mathcal{D}$. Thus, the number of fixed points by $H$ on each $H$-invariant block is $2$. Note that, for any pair of distinct points fixed by $H$ there are exactly $\lambda=3$ blocks of $\mathcal{D}$, and each of these are necessarily left invariant by $H$. So, counting the flags $(\beta,C)$ of $\mathcal{D}$ fixed by $H$ we see that $8 \cdot 2=\binom{4}{2}\cdot \lambda$ with $\lambda=3$, a contradiction. This completes the proof.      
\end{proof}

\begin{lemma}\label{case2}
$X_{\alpha}\neq \mathrm{PSL}(2,q_0)$, for $q=q_0^t$ odd, where $t$ is an odd prime.
\end{lemma}

\begin{proof}
Suppose that $X_{\alpha}=\mathrm{PSL}(2,q_0)$, for $q=q_0^t$ odd, where $t$ is an odd prime.
Let $q_0=p^s$, for some integer $s$, and so $f=st$.
Here, $v=q_0^{t-1}(q_0^{2t}-1)/(q_0^2-1)$ and $|\mathrm{Out}(X)|=2f$.
By $r\mid |G_\alpha|$ and Lemma~\ref{lem:HX}, $r$ divides $fq_0(q_0^2-1)$. 
Since $r^2>3v$ and $q=q_0^t=p^{f}>f^2$, we obtain
\[
    q_0^{t+8} > f^2 q_0^2(q_0^2-1)^3 >3q_0^{t-1}(q_0^{2t}-1) >q_0^{3t-1} ,
\]
implying $t=3$, as $t$ is odd.
Then $f=3s$, $v=q_0^2(q_0^4+q_0^2+1)$ and so $\gcd(v-1,q_0)=1$. Combining this with $r$ divides $\gcd\bigl(3(v-1), fq_0(q_0^2-1)\bigr)$,
we deduce that $r$ divides $3\gcd(v-1, f(q_0^2-1))$,
implying $r\mid 6f$, for $\gcd(v-1, q_0^2-1)=2$.
Therefore, $36f^2 \ge r^2>3v=3q_0^2(q_0^4+q_0^2+1)$, where $f=3s$ and $q_0=p^s$,
and so $108s^2 >p^{6s}$, which is impossible.
\end{proof}

\begin{lemma}\label{case3}
$X_{\alpha}\neq \mathrm{PGL}(2,q_0)$, for $q=2^f=q_0^t$, where $t$ is prime and $q_0\neq 2$.
\end{lemma}

\begin{proof}
Suppose that $X_{\alpha}=\mathrm{PGL}(2,q_0)$, for $q=2^f=q_0^t$, where $t$ is prime and $q_0\neq 2$.
Then $v=q_0^{t-1}(q_0^{2t}-1)/(q_0^{2}-1)$ and $|\mathrm{Out}(X)|=f$, with $q_0=2^s$, for some integer $s$.
Our argument is the same as in the proof of Lemma \ref{case2}. 
By $r\mid |G_\alpha|$ and Lemma~\ref{lem:HX}, $r$ divides $fq_0(q_0^2-1)$. 
Since $r^2>3v$ and $2q=2q_0^t=2^{f+1}>f^2$, we obtain
\[
       q_0^{t+9}>f^2q_0^2(q_0^2-1)^3>3q_0^{t-1}(q_0^{2t}-1)>q_0^{3t-1},
\]
implying $t=2$ or 3.

If $t=2$, then $v=q_0(q_0^{2}+1)$, $q=q_0^{2}$ and $f=2s$. 
By $r\mid \gcd\bigl(3(v-1), fq_0(q_0^{2}-1)\bigr)$ and $\gcd(3(v-1),q_0)=1$, 
we have $r\mid\gcd\bigl(3(v-1), f(q_0^{2}-1)\bigr)$,
and so $r$ divides $3f\gcd\bigl(v-1, q_0^{2}-1)\bigr)$.
As $\gcd(v-1,q_0^{2}-1)=1$ or 3, it follows that $r\mid 9f$.

Thus,
\[
   324s^2=81f^2\ge r^2>3v=3q_0(q_0^{2}+1)>3\cdot 2^{3s},
\]
where $q_0=2^s\ge 4$, implying $s=2$ or 3,
and so $v=68$ with $r\mid 3$, or $v=520$ with $r\mid 9$,
contradicting the fact that $r^2>3v$.

If $t=3$, then $v=q_0^{2}(q_0^{4}+q_0^{2}+1)$, $q=q_0^{3}$ and $f=3s$. 
Again, we have $r$ divides $3f\gcd\bigl(v-1, q_0^{2}-1\bigr)=3f$.
Thus, $27s^2=3f^2>q_0^{2}(q_0^{4}+q_0^{2}+1)>2^{6s}$, which is impossible.

\end{proof}

\begin{lemma}\label{case4}
Suppose that $X_{\alpha}=\mathrm{A_4}$, for $q=p\equiv \pm 3 \text{(mod 8)}$ and $q\not\equiv \pm 1\text{(mod 10)}$. Then $\mathcal{D}$ is a $2$-$(5,3,3)$ design.
\end{lemma}

\begin{proof}
Here, $v=q(q^2-1)/24$ and $|\mathrm{Out}(X)|=2$, and so $r\mid 24$.
By $r^2>3v$, we obtain $8\cdot24^2>q(q^2-1)$,
implying $q=5$, 13, and $v=5$, 91, respectively.

If $(q,v)=(5,5)$, combining with $r\mid 3(v-1)$, $r^2>3v$ and $k\ge 3$,
then $r=6$, $k=3$ and $b=10$, which is a complete design.
If $(q,v)=(13,91)$, then $r^2\le\gcd(270,24)= 6<3v$, a contradiction.

\end{proof}

\begin{lemma}\label{case5}
$X_{\alpha}\neq\mathrm{S_4}$, for $q=p\equiv \pm 1\text{(mod 8)}$.
\end{lemma}

\begin{proof}

Suppose that $X_{\alpha}=\mathrm{S_4}$, for $q=p\equiv \pm 1\text{(mod 8)}$.
Then $v=q(q^2-1)/48$ and $|\mathrm{Out}(X)|=2$, and so $r\mid 48$.
Again, by $r^2>3v$, we obtain $16\cdot48^2>q(q^2-1)$,
implying $q=7$, 17, 23, 31, and so $v=7$, 102, 253, 620, respectively. 
If $v=7$, then by $r\mid\gcd(48,18)$ and $r^2>3v$, we have 
$r=6$, and so $k=4$, contradicting the fact that $b$ is an integer. 
For the remaining cases, 
we always have $r^2\le\gcd^2(48,3(v-1))\le 3v$, a contradiction.

\end{proof}

\begin{lemma}\label{case6}
Suppose that $X_{\alpha}=\mathrm{A_5}$, for $q\equiv \pm 1\text{(mod 10)}$, Where either $q=p$ or $q=p^2$ and $p\equiv \pm 3\text{(mod 10)}$. Then $G=\PSL(2,11)$ and $\mathcal{D}$ is isomorphic to the $2$-$(11,3,3)$ as in Example \ref{Ex1}.
\end{lemma}

\begin{proof}
Then $v = q(q^2-1)/120$, $|\mathrm{Out}(X)|=2f$, and so $r\mid 120f$, with $f=1$ or 2.
By $r^2>3v$, we obtain $40 \cdot120^2f^2>q(q^2-1)$.

Consider first $f=2$, then $q=p^2\le 132$, and so $(q,v)=(9,6)$ or (49,980). 
If $(q,v)=(49,980)$, then by $r\mid \gcd(240, 3\times979)$, we have $r\mid3$, a contradiction.
If $(q,v)=(9,6)$, then $r\mid 15$. Since $r^2>3v$ and $r>k\ge 3$,  we have $r=5$ and $k=4$, contradicting the fact that $b$ is an integer.
Thus, $f=1$.  Then $q= p\le 83$, by $v = q(q^2-1)/120$ and $q \equiv \pm 1 \pmod{10}$, we obtain
$(q,v)=(11,11)$, (19,57), (29,203), (31,248), (41,574), (59,1711),
(61,1891), (71,2982) or (79,4108).
Combining this with  $r\mid \gcd(120,3(v-1))$ and $r^2>3v$,
we have that $(q,v)=(11,11)$ or $(19,57)$,
and so $(v,b,r,k)=(11,55,15,3)$ or $(57,171,24,8)$, respectively. 

In the latter case, the flag-transitivity implies that the order of $G$ is divisible by $vr=57\cdot 24$, forcing $G=\PGL(2,19)$. Then $\left \vert G_{\alpha}\right\vert =120$ containing $X_{\alpha}=\mathrm{A_5}$. However, $G$ does not have such subgroups. Indeed, $\PSL(2,19)$ has two conjugate classes of (maximal) subgroups isomorphic to $\mathrm{A_5}$, and these are fused in $\PGL(2,19)$.  

Finally, if $(v,b,r,k)=(11,55,15,3)$ then $G=X$ acts point-$2$-transitively on $\mathcal{D}$, and hence either $X_{B} \cong \mathrm{A}_{4}$ or $\mathrm{D}_{12}$ by \cite{At}. In the former case, $X_{B} \leq X_{\alpha} \cong \mathrm{A}_{5}$.  for some $\alpha \notin B$ by \cite{At}. 
Hence, $\mathrm{A_5}$ acts on the $10$ points of $\mathcal{D}$ distinct from $\alpha$. 
Now, $\mathrm{A_5}$ has exactly one transitive permutation representation of degree 10 by \cite{At}, and this one is equivalent to the action of $\mathrm{A_5}$ in $\PG(1,9)$ when regarded as subgroup of $\PSL(2,9)$ by \cite[Lemma 11.(ii)]{COT}. Since $\mathrm{A_4}$ partitions $\PG(1,9)$ into two orbits of length $6$ and $4$, $X_{B} \cong \mathrm{A_4}$ does not have orbits of length $3$ on the set of points of $\mathcal{D}$ distinct from $\alpha$. Therefore, $X_{B} \cong D_{12}$.

Let $C$ be the normal subgroup of $X_{B}$ of order $3$. Then $C$ is a Sylow $3$-subgroup of $X$, and $N_{X}(C)=X_{B}$ since $X_{B}$ is a maximal subgroup of $G$. Now, $C$ fixes at least one point as $v=11$, that we may assume to be $\alpha$. Since $N_{X_{\alpha}}(C) \cong S_{3}$, it follows that $C$ fixes exactly $\left\vert N_{X}(C): N_{X_{\alpha}}(C) \right\vert=2$ points of $\mathcal{D}$, and these are permuted by $X_{B}$. 
The remaining $9$ points of $\mathcal{D}$ are partitioned into exactly $3$ orbits, each of length $3$. 
At this point, it is not difficult to see that $X_{B}$ preserves exactly one $C$-orbit of length $3$ and switches the remaining two, since an elementary abelian group of order $4$ does not fix $3$ or more points. Thus, $\mathcal{D}$ is isomorphic to the $2$-$(11,3,3)$ design as in Example \ref{Ex1}.   
\end{proof}

\begin{lemma}\label{case7}
$X_{\alpha}\neq\mathrm{D}_{2(q-1)/(2,q-1)}$.
\end{lemma}

\begin{proof}
Suppose that $X_{\alpha}=\mathrm{D}_{2(q-1)/(2,q-1)}$. 
Then $v=\frac{q(q+1)}{2}$. 

If $q$ is even, by Lemma \ref{subdegrees-dihedral}, then we obtain that $r$ divides $\frac{3(q-2)(q-1)}{2}$.
Combining this with $r\mid 3(v-1)$, we get $r\mid 3(q-1)$.
It follows from $r^2>3v$ that $r=3(q-1)$, and so $k=\frac{q+4}{2}$, $b=\frac{3q(q^2-1)}{q+4}$.
As $b$ is an integer, we have that $q+4$ divides $180$, implying $q=8$, $16$ or $32$, 
and so
\[
  (v,b,r,k)=(36,126,21,6), (136,612,45,10) \textrm{~or~}(528,2728,93,18).
\]
The first case is ruled out by \cite[Theorem 1.1]{MF}. 
For the remaining cases, $G \leq \mathrm{P \Gamma L}(2,16)$ or $G \leq \mathrm{P \Gamma L}(2,32)$, respectively. 
Since $G$ is flag-transitive, we have that $v \cdot r$ divides $|G|$, and so $9\mid |G|$, a contradiction.
Thus, $q$ is odd. Similarly, by Lemma \ref{subdegrees-dihedral} and $r^2>3v$, we obtain $r=\frac{3(q-1)}{2}$,
and so $k=q+3$, $b=\frac{3q(q^2-1)}{4(q+3)}$,  which implies $q+3\mid 72$.
Hence, $q=5$ or 9, contradicting the fact that $r>k$.
\end{proof}

\bigskip
Before proceeding, we need to recall the following facts related to $\PSL(2,q)$, $q$ even, as a collineation group of $\PG(2,q)$:
\begin{enumerate}
\item[(1)] An irreducible conic $\mathcal{C}$ of $\PG(2,q)$, $q=2^{f}$, is a $%
(q+1)$\emph{-arc}, namely a set of $q+1$ points no three of them collinear, by 
\cite[Lemma 7.7]{Hir}. Any line of $\PG(2,q)$ is either \emph{secant}, \emph{%
tangent} or \emph{external} according as it has $2$, $1$ or $0$ points in
common with $\mathcal{C}$, respectively. The set of secant, tangent or external lines has size $\frac{q(q+1)}{2}$, $q+1$, or $\frac{q(q-1)}{2}$, respectively, by \cite[Corollary 8.2]{Hir}. The tangent lines to $\mathcal{C}$
are all concurrent to a point $N$ called \emph{nucleus} of $\mathcal{C}$ by 
\cite[Corollary 7.11]{Hir}, and $\mathcal{J}=\mathcal{C}\cup \left\{
N\right\} $ is a $(q+2)$-arc called \emph{regular} \emph{hyperoval }(see \cite[%
Section 8.4]{Hir}). Clearly, the lines of $\PG(2,q)$ are either secants or external to $\mathcal{J}$, that is they have either $2$ or $0$ points in common with $%
\mathcal{J}$. The set $\mathcal{E}$ of the external lines to $\mathcal{J}$ coincides with the set of the external lines to $\mathcal{C}$, and hence it has size $q(q-1)/2$. The number of
points of $\PG(2,q)\setminus \mathcal{J}$ is $q^{2}-1$ and through each point
there are exactly $\frac{q}{2}+1$ secant lines to $\mathcal{J}$ and $q/2$ external lines to $\mathcal{J}$ by \cite[Corollary 8.8]{Hir}.

\item[(2)] $\PGL(3,q)$ has a unique conjugacy class of subgroups isomorphic
to $\PSL(2,q)$ by \cite[Table 8.3]{BHRD} (note that $\PSL(2,q) \cong \Omega (3,q)$ by \cite[Corollary 7.14]{Hir} is reducible
and not maximal in $\PGL(3,q)$ when $q$ is even. Further, the irreducible conics do not arise from
polarities in this case), and each of these groups is the stabilizer in $%
\PGL(3,q)$ of a suitable regular hyperoval of $\PG(2,q)$. The converse is also
true as a consequence of \cite[Theorem 7.4]{Hir}. In particular, each $%
\PSL(2,q)$ fixes the nucleus of its invariant hyperoval and acts $2$%
-transitively on the remaining $q+1$ points of this one  by \cite[Corollary 7.15]{Hir}.

\item[(3)] $\PGL(3,q)$ has a unique conjugacy class of elements of order $2$, and if $%
\sigma $ is any of these then $\sigma $ is a $(P_{\sigma },t_{\sigma })$%
-elation of $\PG(2,q)$. That is, $\sigma $ fixes each of the $q+1$ points of $%
t_{\sigma }$ including $P_{\sigma }$ and fixes setwise each of the $q+1$ lines
of $\PG(2,q)$ containing $P_{\sigma }$ by \cite[Exercise IV.4.6]{HP}. No
other points or lines of $\PG(2,q)$ are fixed by $\sigma $.
\item[(4)] If $T$ is any Sylow $2$-subgroup $S$ of $\PSL(2,q)$, then the following hold by \cite[Lemma 2.3]{MZZZ}:
\begin{enumerate}
    \item $T$ fixes the nucleus $N$, a unique point $Q$ of $\mathcal{C}$ and acts regularly on $\mathcal{C}\setminus \left\{ Q\right\}$;
    \item $T$ fixes $t$ pointwise, where $t=NQ$, and acts semiregularly on $\PG(2,q)\setminus t$;
    \item  $\sigma _{i}$ is a $(P_{\sigma_{i}},t)$-elation of $\PG(2,q)$;
    \item  Any cyclic subgroup of $\PSL(2,q)$ normalizing $T$ acts regularly on the set $t\setminus \left\{ N,Q\right\}=\left\{ P_{\sigma _{1}},...,P_{\sigma _{q-1}}\right\}$;
    \item If $E_{i}$ is the set of $q/2$ lines through $P_{\sigma_{i}}$ which are external to $\mathcal{J}$ (see (2)), then $T$ acts transitively on $E_{i}$ with action kernel $\left\langle \sigma _{i} \right\rangle$. Moreover, $\mathcal{E}$ is partitioned into $q-1$ orbits under $T$, namely, $E_{1},...,E_{q-1}$, which are the unique $S$-orbits.
\end{enumerate}
\end{enumerate}

\bigskip

\begin{lemma}\label{case8}
$X_{\alpha}\neq\mathrm{D}_{2(q+1)/(2,q-1)}$.
\end{lemma}

\begin{proof}
Suppose that $X_{\alpha}=\mathrm{D}_{2(q+1)/(2,q-1)}$. 
Then $v=\frac{q(q-1)}{2}$. 
Let $B$ be a block of $\mathcal{D}$ incident with point $\alpha$.
Since $G$ is block transitive and $X\triangleleft G$, we have that 
$|X:X_{B}| \mid b$. We may assume that $N$ is a maximal subgroup of $X$ containing $X_B$,
and so $|X:N|\mid b$. We now analyze all the possibilities of $N$.

First, consider the case that $q$ is odd. Then by Lemma \ref{subdegrees-dihedral}, we have $r$ divides $3(q+1)/2$.
It follows from $r^2>3v$ that $r=3(q+1)/2$, and so $k=q-1$, $b=\frac{3q(q+1)}{4}$,
implying $q\equiv 3\pmod{4}$. If $N=\mathrm{PGL}(2,q_0)$, then $q=q_0^2 \equiv 3\pmod{4}$, which is impossible.

Suppose that $N=\PSL(2, q_0)$ with $q=q_0^t$ odd, where $t$ is odd and prime.
Then $|X:N|=q(q^2-1)/q_0(q_0^2-1)$, and so $4(q-1) \mid 3q_0(q_0^2-1)$.
Thus, $4(q_0^3-1)\le 4(q-1) \le 3q_0(q_0^2-1)$, a contradiction.

Suppose that $N=\mathrm{A}_4$, for $q=p\equiv\pm3\pmod8$ and $q\not\equiv\pm1\pmod{10}$.
Then $q(q^2-1)/24$ divides $3q(q+1)/4$, and so $q-1\mid 18$, implying $q=7$ or 19, a contradiction. 

Suppose that $N=\mathrm{S}_4$, for $q=p\equiv\pm1\pmod8$.
Then $q(q^2-1)/48$ divides $3q(q+1)/4$, and so $q-1\mid 36$, forcing $q=7$, 
and so $(v,b,r,k)=(21,42,12,6)$. 
It follows from $\PSL(2,7)\unlhd G \leq \PGL(2,7)$ and $v=21$ that $G_{\alpha}$ is a Sylow $2$-subgroup of $G$, contradicting the fact that $r=12$ divides $|G_{\alpha}|$.

Suppose that $N=\mathrm{A}_5$, for $q\equiv \pm 1\text{(mod 10)}$, where either $q=p$ or $q=p^2$ and $p\equiv \pm 3\pmod{10}$.
Then $q(q^2-1)/120$ divides $3q(q+1)/4$, and so $q-1\mid 90$,
implying that $q=11$, $19$ or $31$. Then $G=\PSL(2,q):Z_{2^{\epsilon}}$ with $\epsilon=0,1$, and hence 
\[
  (v,b,r,k)=(55,99,18,10), (171,285,30,18) \textrm{~or~}(465,744,48,30).
\]
Then $\left \vert G_{\alpha} \right\vert =12 \cdot 2^{\epsilon}, 20 \cdot 2^{\epsilon}$, or $32 \cdot 2^{\epsilon}$, respectively, and these cases are ruled out since they are not divisible by $r=18, 30$ or $48$, respectively.

If $N=\mathrm{D}_{q-1}$, then $q(q+1)/2$ divides $3q(q+1)/4$, a contradiction. 
If $N=\mathrm{D}_{q+1}$ then $q(q-1)/2$ divides $3q(q+1)/4$, and so $q-1\mid 6$, forcing $q=7$. Hence, $(v,b,r,k)=(21,42,12,6)$ and $\PSL(2,7) \unlhd G \leq \PGL(2,7)$. The flag-transitivity of $G$ implies that the order of $G$ must be divisible by $v\cdot r= 2^{2}\cdot 3^{2} \cdot 7$, which is not the case, and hence there is no such design. If $N=C_p^f\rtimes C_{q-1}/2$, then $q+1$ divides $3q(q+1)/4$, which is impossible.

Therefore, $X_{\alpha}=\mathrm{D}_{2(q+1)}$ and $|\Out(X)|=f$ with $q=2^f$. 
If $q=4$, then $v=6$. However, there is no $r$ satisfying $r\mid 3(v-1)$, $r>k>2$ and $3\mid r$. Similarly, if $q=8$, then $(v,b,r,k)=(28,189,27,4)$. Moreover, $G=\mathrm{P \Gamma L}(2,8)$ since $r=27$ divides the order of $G$. Now, since $G$ acts flag-transitively on $\mathcal{D}$, we see that $G_{x,B}=Z_{2}$. On the other hand, the order of $G_{B}$ is divisible by $4$ since $k=4$, and hence $G_{B}$ contains a Sylow $2$-subgroup, say $S$, containing $G_{x,B}=Z_{2}$. Since $S$ is elementary abelian of order $8$ and acts transitively on $B$, it follows that $G_{x,B}=Z_{2}$ fixes $B$ pointwise. The action on the point set of $\mathcal{D}$ and the one on $\mathcal{R}(3)$, the Ree unital of order $3$, are equivalent, and hence the fixed points of an involution are exactly a block of $\mathcal{R}(3)$ (see \cite{Lu}), forcing $\mathcal{D} \cong \mathcal{R}(3)$ since $G$ is block-transitive on both designs. This is a contradiction since $\mathcal{R}(3)$ is linear space, whereas $\lambda=3$. So, there is no such design $\mathcal{D}$.
Hence, we may assume that $q=2^f\ge 16$.
By Lemma \ref{lem:HX} and Lemma \ref{subdegrees-dihedral}, we deduce 
\[
r\mid (q+1)\gcd(3(2^{f-1}-1),2f).
\]
Let $r=\frac{n(q+1)}{m}$, where $m$ and $n$ are positive integers such that
$\gcd(m,n)=1$, $m\mid(q+1)$ and  $n\mid \gcd(3(2^{f-1}-1),2f)$.
Then $m$, $n$ are odd, and 
\[
k=\frac{3m(q-2)+2n}{2n} \textrm{~and~} b=\frac{n^2q(q^2-1)}{m\big(3m(q-2)+2n\big)}.
\]
Since $k<r$, $3m^2(q-2)<2n^2(q+1)$, we conclude $m<n$, as $q\ge16$.
From the fact that $|X:N|$ divides $b$,  we analyze all the possibilities of $N$.

\noindent\textbf{(9.1)} If $N=\mathrm{D}_{2(q-1)}$, then $q(q+1)/2$ divides $n^2q(q^2-1)/[m(3m(q-2)+2n)]$,
and so $m(3m(q-2)+2n)$ divides $2n^2(q-1)$. Since $\gcd(m,n)=1$ and $m$ odd, we have $m\mid (q-1)$,
and so $m\mid \gcd(q-1,q+1)=1$, $m=1$, $n\ge3$.
Thus, $r=n(q+1)$, $k=[3(q-2)+2n]/2n$ and $b=n^2q(q^2-1)/[3(q-2)+2n]$, and so
\[
3(q-2)+2n \mid 2n^2(q-1).
\]
Set $E:= 3(q-2)+2n$. Then $E \mid 6n^{2}(q-1)$. 
As $6n^{2}(q-1)=2n^{2}E+2n^{2}(3-2n)$, we have $E\mid 2n^{2}(2n-3)$.
Thus, $3(q-2)<E\le2n^{2}(2n-3)<4n^{3}$, and so
\[
2^{f}=q<2n^{3}\le 2f^{3},
\]
as $n\mid \gcd\bigl(3(2^{f-1}-1),2f\bigr)$ and $n$ is odd.
Hence, $4\le f\le 11$. 
Combining with the fact that odd $n$ divides $\gcd\big(3(2^{f-1}-1),\,2f\big)$ 
and $2^{f}<2n^{3}$, we obtain
\[
(f,n)\in\{(5,5),(7,7),(9,9),(11,11)\}.
\]
However, none of them satisfy $E\mid 2n^{2}(q-1)$.

\noindent\textbf{(9.2)} If $N=\mathrm{D}_{2(q+1)}$, then $q(q-1)/2$ divides $n^2q(q^2-1)/[m(3m(q-2)+2n)]$,
and so $m(3m(q-2)+2n)$ divides $2n^2(q+1)$.
Set $E:=3m(q-2)+2n$, and so $E\mid 2n^2(q+1)$.

Assume that $m=1$. Then 
\[
E=3(q-2)+2n \mid 2n^2(q+1).
\]
If $n=3$, then $3q\mid 18(q+1)$, a contradiction.
Hence, $n\ge 5$.
Note that $6n^2(q+1)=2n^2E+2n^2(9-2n)$. We deduce that $E\mid 2n^2(2n-9)$, 
and so $3(q-2)<6n^3$, implying 
\[
2^{f}=q<3n^{3}\le 3f^{3}.
\]
Therefore, $4\le f\le 12$. 
Similarly, by $n$ divides $\gcd\big(3(2^{f-1}-1),\,2f\big)$ 
and $2^{f}<3n^{3}$, we obtain
\[
(f,n)\in\{(5,5),(7,7),(9,9),(11,11)\}.
\]
Again, none of them satisfy $E\mid 2n^{2}(q+1)$.

Thus, $m\ge3$ and $n\ge5$.
Since $6mn^{2}(q+1)= 2n^{2}E+2n^{2}(9m-2n)$, 
we have that $E$ divides $2n^{2}(9m-2n)$, and so 
$9(q-2)<E\le 2n^{2}|9m-2n|<14n^{3}$,
implying 
\[
2^{f}=q<2n^{3}\le2f^{3}.
\]
We get 
\[
(f,n)\in\{(5,5),(7,7),(9,9),(11,11)\}.
\]
However, there is no odd $m$ satisfying $m(3m(q-2)+2n)\mid [2n^2(q+1)]$ with $\gcd(m,n)=1$ 
and $3\le m<n$.

\noindent\textbf{(9.3)} Suppose that $N=\PGL(2,q_{0})$ with $q=2^{f}=q_{0}^{t}$ for a prime $t$ and $q_0 \neq 2$. Then $f\ge 2t$, $q_0\ge4$.
Moreover, from the fact that $|X:N|$ divides $b$, 
we have 
\[
m\bigl(3m(q-2)+2n\bigr)\mid n^{2}q_{0}(q_{0}^{2}-1),
\]
and so $3m^2q\le n^{2}q_{0}(q_{0}^{2}-1)+\lvert 6m^{2}-2mn\rvert$.
It follows from $\lvert 6m^{2}-2mn\rvert<4n^2$ that $3m^2q_0^t <n^2q_0^3$.
If $t\ge 5$, then $3m^2q_0^{t-3} <n^2$.
Note that $2^{f}=q_{0}^{t}$, $n$ odd and $n\mid \gcd\bigl(3(2^{f-1}-1),2f\bigr)$.
We deduce 
\[
 2^{2f/5}<2^{f(t-3)/t}=q_0^{t-3} <n^2\le f^2,
\]
and so $2^f<f^5$, implying $10\le f\le 22$, for $f\ge 2t$.
Combining this with prime $t\ge5$ and $q_0\ge 4$, we have
$f=10$, 14, 15, 20, 21, or 22, and so $n\leq \gcd\bigl(3(2^{f-1}-1),2f\bigr)\leq 3$,
contradicting the fact that $n^2>3m^2q_0^{t-3}\ge48$.
Therefore, $t=2$ or 3. 

If $t=2$, then $m\bigl(3m(q_{0}^{2}-2)+2n\bigr)\mid n^{2}q_{0}(q_{0}^{2}-1)$.
Since $m\mid q_0^2+1$ and $\gcd(m,n)=1$, we have $m\mid q_{0}^{2}-1$, 
and so $m\mid \gcd(q_{0}^{2}-1,q_{0}^{2}+1)$, forcing $m=1$.
Hence, $n\ge 3$ and $3(q_{0}^{2}-2)+2n \mid n^{2}q_{0}(q_{0}^{2}-1)$.
If $n=3$, then $3q_{0}^{2}\mid 9q_{0}(q_{0}^{2}-1)$, which is impossible.
So $n\ge 5$. Set $E:=3(q_{0}^{2}-2)+2n$.
Note that $3n^{2}q_{0}(q_{0}^{2}-1)=n^{2}q_{0}E + n^{2}q_{0}(3-2n)$.
We deduce $E\mid n^{2}q_{0}(2n-3)$, and so $3q_{0}^{2}<E<2n^3q_0$,
implying 
\[
2^{f/2}=q_{0}<n^{3}\le f^3.
\]
Therefore, $2^{f}<f^{6}$, $4\le f\le 28$ with $f$ even.
Combining this with $n^3\ge 2^{f/2}$ and $n\mid\gcd\bigl(3(2^{f-1}-1),\,2f\bigr)$,
we have $(f,n)=(6,3)$, contradicting the fact that $n\ge 5$.

If $t=3$, then $m\bigl(3m(q_{0}^{3}-2)+2n\bigr)\mid n^{2}q_{0}(q_{0}^{2}-1)$.
Set $E:=m\bigl(3m(q_{0}^{3}-2)+2n\bigr)$.
Similarly, from 
$3m^{2}n^{2}q_{0}(q_{0}^{2}-1)=n^{2}E-mn^{2}\bigl(3mq_{0}-6m+2n\bigr)$, 
we obtain that $E\mid mn^{2}\bigl(3mq_{0}-6m+2n\bigr)$, then 
$2mq_{0}^{3}<3mn^{2}q_{0}+2n^{3}$,
and so 
\[
2\cdot 2^{2f/3}=2q_{0}^{2}<3n^{2}+\frac{2n^{3}}{mq_0}\le n^3+\frac{n^{3}}{2}\le \frac{3f^{3}}{2},
\]
implying $4^f<f^9$.
Therefore, $f=6$, 9, 12, 15 or 18.
Recall that $\frac{3n^{3}}{2}>2\cdot 2^{2f/3}$ and $n\mid\gcd\bigl(3(2^{f-1}-1),\,2f\bigr)$,
we have $(f,n)=(6,3)$ or (9,9), and so $m=1$, as $m<n$, $m\mid q+1$ and $\gcd(m,n)=1$.
However, neither of them satisfy $m\bigl(3m(q_0^3-2)+2n\bigr)$ divides $n^2 q_0(q_0^2-1)$.

\noindent\textbf{(9.4)} Suppose that $N= q:(q-1)$.  
Then $m\big(3m(q-2)+2n\big)$ divides $n^2q(q-1)$,
and so $m\mid \gcd(q-1,q+1)$, implying $m=1$.
Thus, $n\ge 3$ and $3(q-2)+2n\mid n^2q(q-1)$.
Set $E:=3(q-2)+2n$.
Note that $9n^{2}q(q-1)=(-2n^{3}+3n^{2}q+3n^{2})E+(4n^{4}-18n^{3}+18n^{2})$.
If $n\ge5$, then $E\mid 2n^{2}(2n^{2}-9n+9)$,
and so $2q<E\le 2n^{2}(2n^{2}-9n+9)<4n^4$, 
implying
\[
2^f=q<2n^4\le 2f^4.
\]
Therefore, $5\le f \le 17$.
Combining with the fact that odd $n\mid\gcd\big(3(2^{f-1}-1),\,2f\big)$,
$2^{f}<2n^{4}$ and $E\mid n^{2}(q-1)$, we obtain $(f,n)=(5,5)$ or (7,7),
and so $q=2^5, 2^7$ 
\[
(v,b,r,k)=(496,8184,165,10)\textrm{~or~}(8128,262128,903,28)
\]
respectively. Hence, either $G = \mathrm{P\Gamma L}(2,2^5)$ with $\left\vert G_{B}\right\vert=10$, or $G=\mathrm{P\Gamma L}(2,2^7)$ with $\left\vert G_{B}\right\vert=56$. If $G_{B}$ contains a normal subgroup $J$ of order $u$, where $u=5$ or $7$, then $N_{G}(J)$ has order divisible by $4$ or $8$, respectively. However, this is impossible since $J$ fixes a subline $\PG(1,2)$ of the projective line $\PG(1,q)$, $q=2^5,2^7$, and hence $N_{G}(J)=J \times S_{3}$ since the group induced by the stabilizer of $G$ on any copy of $\PG(1,2)$ is $\PGL(2,2)\cong S_{3}$ and the stabilizer in $\PGL(2,q)$ of any three distinct points of $\PG(1,q)$ is trivial. Therefore, $G_{B}$ does not contain normal subgroups of order $u$, where $u=5$ or $7$. Then $G=\mathrm{P\Gamma L}(2,2^7)$, $G_{\alpha}=D_{258}:C_{7}$, and $G_{B}=(C_{2})^{3}:C_{7}$. Then $G_{\alpha,B}=C_{2}$ fixes exactly $4$ points in $B$. 

By (3), we may assume that $G$ is the copy of $\mathrm{P\Gamma L}(2,2^7)$ inside 
$\mathrm{P \Gamma L}(3,2^{7})$
preserving a fixed hyperoval $\mathcal{J}$ and, as pointed out in \cite[Section 2.6]{BDD}, we may identify
the point set of $\mathcal{D}$ with $\mathcal{E}$, the set of lines $\PG(2,2^{7})$ external to $\mathcal{J}$. Now any element of order $7$ fixes a Fano subplane $ \pi =\PG(2,2)$ of $\PG(2,2^7)$. Easy computations show that $\pi$ has exactly $4$ points in $\mathcal{J}$ and $6$ lines of $\PG(2,2^7)$ secants to $\mathcal{J}$. Thus, any element of order $7$ fixes exactly one line external to $\mathcal{J}$, that is, any element of order $7$ fixes exactly one point of $\mathcal{D}$. This, together with $X_{\alpha}$ partitions the points of $\mathcal{D}$ distinct from $\alpha$ into $63$ orbits of length $129$ by Line 2 in Table \ref{tab:subdegrees-PSL2-dihedral}, implies that $G_{\alpha}$ partitions the points of $\mathcal{D}$ distinct from $\alpha$ into $9$ orbits of length $903$. Let $\mathcal{O}_{i}$, $i=1,...,9$ any such $G_{\alpha}$-orbit. Then $r\left\vert B \cap \mathcal{O}_{i}\right\vert=3\left\vert \mathcal{O}_{i}\right\vert$ by \cite[(1.2.6)]{DembowskiFG} since $\lambda=3$. Now, since $r=\left\vert \mathcal{O}_{i}\right\vert=903$ for each $i=1,...,9$, it follows that $\left\vert B \cap \mathcal{O}_{i}\right\vert=3$ for each $i=1,...,9$. Then $G_{\alpha,B}=C_{2}$ fixes at least a point on $B \cap \mathcal{O}_{i}$ for each $i=1,...,9$, since each $B \cap \mathcal{O}_{i}$ is left invariant by $G_{\alpha,B}$. So, $G_{\alpha,B}$ fixes at least $10$ points $B$, including $\alpha$. However, this is impossible since we have seen that $G_{\alpha,B}=C_{2}$ fixes exactly $4$ points in $B$. Thus, this case is ruled out.   

Therefore, $m=1$ and $n=3$. We have 
\[
(v,b,r,k)
=\bigl(q(q-1)/2,\;3(q^{2}-1),\;3(q+1),\;q/2\bigr).
\]
Since $X_{\alpha }=\mathrm{D}_{2(q+1)}$ and $\left\vert X:X_{B}\right\vert =\frac{%
3(q^{2}-1)}{\theta }$ and so $X_{B}$ contains a Sylow Subgroup of $X$.
Therefore, $X_{B}=q:\theta $ with $\theta \mid 3(q-1)$. Since $X$ contains a
unique conjugacy class of subgroups isomorphic to $q:\theta $, then $%
\left\vert G_{B}/X_{B}\right\vert =\left\vert G/X\right\vert $, hence 
\[
3(q^{2}-1)=b=\left\vert G:G_{B}\right\vert =\left\vert X:X_{B}\right\vert =%
\frac{3(q^{2}-1)}{\theta }
\]%
forces $\theta =1$. So, $X_{B}$ is a Sylow $2$-subgroup of $X$ and $X$ acts
block-transitively on $\mathcal{D}$. Moreover, 
\[
q/2=k\geq \left\vert X_{B}:X_{\alpha ,B}\right\vert \geq q/2
\]%
since $X_{\alpha }=\mathrm{D}_{2(q+1)}$ and therefore, $X$ acts flag-transitively
on $\mathcal{D}$. Thus, we may assume that $G=X$ without loss.

By (3), in the sequel, we may assume that $X$ is the copy of $\PSL(2,q)$ inside $\PGL(3,q)$
preserving a fixed $\mathcal{J}$ and, as pointed out in \cite[Section 2.6]{BDD}, we may identify
the point set of $\mathcal{D}$ with $\mathcal{E}$, the set of lines 
$\PG(2,q)$ external to $\mathcal{J}$. Moreover, if $B$ is any block of $\mathcal{D}$ then $X_{B}$ is a Sylow $2$-subgroup, and it follows from (5) that $B$ is a pencil of $q/2$ lines of $\PG(2,q)$ external to $\mathcal{J}$. Now, let $\ell$ and $\ell^{\prime}$ be any two distinct points of $\mathcal{D}$. Then $\ell$ and $\ell^{\prime}$ are elements of $\mathcal{E}$. Therefore, $\ell$ and $\ell^{\prime}$ are two distinct lines of $\PG(2,q)$ external to $\mathcal{J}$.

If there are three distinct blocks of $\mathcal{D}$ containing both $\ell$ and $\ell^{\prime}$, say $B_{1},B_{2},B_{3}$, then these are three distinct pencils of $q/2$ lines of $\PG(2,q)$ external to $\mathcal{J}$. Then the centers of the pencils of $B_{1},B_{2},B_{3}$, say $P_{1},P_{2},P_{3}$, respectively, belong to $\ell \cap \ell^{\prime}$. Therefore, either $P_{1}=P_{2}=P_{3}$ or $\ell=\ell^{\prime}$. The latter is clearly ruled out since $\ell$ and $\ell^{\prime}$ are distinct. Therefore, $P_{1}=P_{2}=P_{3}$ and hence $B_{1}=B_{2}=B_{3}$ by (2), a contradiction. 

\end{proof}

\begin{lemma}\label{case9} Suppose that $X_{\alpha}=C_p^f\rtimes C_{q-1/\gcd(2,q-1)}$.
Then one of the following holds:
\begin{enumerate}
    \item $\mathcal{D}$ is the complete $2$-$(5,3,3)$ design 
    and $\PSL(2,4)\unlhd G \leq \mathrm{P \Gamma L}(2,4)$;
    \item $\mathcal{D}$ is the $2$-$(26,6,3)$ design as in Example \ref{Ex2} and 
    $\PSL(2,25)\unlhd G\leq \mathrm{P\Sigma L}(2,25)$.  
\end{enumerate}
\end{lemma}

\begin{proof} 
If $X_{\alpha}=C_p^f\rtimes C_{q-1/\gcd(2,q-1)}$, then the stabilizer $X_{\alpha,\beta}$ has orbit sizes 1, 1, $(q-1)/2$, $(q-1)/2$ if $q$ odd, 1, 1, $q-1$ if $q$ even. 
Here $v=q+1$, and so $r\mid 3q$. 
As $3\leq k<r$ and $3\mid r$, we may assume $r=3p^x$ for some positive integer $x$. 
It follows that $k=p^{f-x}+1$ with $f\ge x+1$. 
Since $r^2>3v$, we obtain $p^{2x+2}> 3p^{2x}>p^f$, and so $f\le 2x+1$. 
Therefore, $2\le x+1\leq f\leq 2x+1$. From $bk=vr$, we have 
\[b=\frac{3p^x(q+1)}{p^{f-x}+1}.\]
Let $B$ be a block of $\mathcal{D}$ incident with point $\alpha$.
Similarly, we may assume that $N$ is a maximal subgroup of $X$ containing $X_B$,
and so $|X: N|$ divides $b$. We now analyze all the possibilities of $N$.

Since $f\ge 2$, we have $N\ne \mathrm{A}_4$ and $\mathrm{S}_4$.
If $N=\mathrm{A}_5$, then $p\equiv\pm3\pmod {10}$, $f=2$ and $x=1$, implying $r=3p$ and $b=3p(p^2+1)/(p+1)$. 
As $b$ is an integer, we have $p+1$ divides $3(p^2+1)$, forcing $(p+1)\mid 6$, 
and so $p=2$ or 5, a contradiction.

Suppose that $N=\PGL(2,q_0)$ with $q=q_0^2$ odd.
Let $q_0=p^s$, and so $f=2s$.
Since $|X:N|$ divides $b$, we have 
\[\frac{q_0(q_0^2+1)}{2}\mid \frac{3p^x(q_0^2+1)}{p^{f-x}+1},\]
and so $p^{f-x}+1\mid 6$, implying $p=5$, $f=x+1$ and $k=6$.
Let $\alpha$ and $\beta$ be two distinct points of $\mathcal D$. 
Since $\lambda=3$, there are exactly three blocks $B_1,B_2,B_3$ containing $\alpha$ and $\beta$.
Moreover, $G_{\alpha,\beta}$ fixes $B_1\cup B_2\cup B_3$ setwise, 
and so $B_1\cup B_2\cup B_3$ is a union of $G_{\alpha,\beta}$-orbits. 
Since $X_{\alpha,\beta}$ has orbit sizes 1, 1, $(q-1)/2$, $(q-1)/2$,
we have that the orbit length of $G_{\alpha,\beta}$ on
$\mathcal P\setminus\{\alpha,\beta\}$ is at least $(q-1)/2$.
Hence, $|B_1\cup B_2\cup B_3|\;\ge\;2+\frac{q-1}{2}$.
Combining this with the fact that $|B_1\cup B_2\cup B_3|\le 14$,
we conclude that $q\le 25$, implying $f=2$. Then $\mathcal D$ has parameters $(v,b,r,k)=(26,65,15,6)$. Moreover $G$ acts point $2$-transitively on $\mathcal{D}$, and the stabilizer in $X$ of a block $B$ of $\mathcal{D}$ is $\PGL(2,5)$. Then we may identify the point set of $\mathcal{D}$ with that of $\PG(1,25)$, now any $\PGL(2,5)$ partitions $\PG(1,25)$ into two orbits, one of them being a Baer subline $\PG(1,5)$, which has size $6$, and its complementary set (of size $20$). Thus $B$ is a copy of $\PG(1,5)$, and hence $\mathcal{D}$ is the $2$-$(26,6,3)$ design as in Example \ref{Ex2}, and so $\PSL(2,25) \unlhd G \leq \mathrm{P\Sigma L}(2,25)$. 

Suppose that $N=\PGL(2,q_0)$ with $q=2^f=q_0^t$,  where $t$ is prime and $q_0\neq 2$;
or $N=\PSL(2,q_0)$ with $q=q_0^t$ odd, where $t$ is odd and prime.
Let $q_0=p^s$, and so $f=st$. Here we have 
\[\frac{q_0^{t-1}(q^{2}-1)}{q_0^2-1}\mid \frac{3p^x(q+1)}{p^{f-x}+1},\]
implying $(p^{f-x}+1)(q-1)$ divides $3(q_0^2-1)$.
Then 
\[
p^{2st-x}<(p^{st-x}+1)(p^{st}-1)<3p^{2s}<p^{2s+2},
\]
and so $2st-2s<2+x\le st+1$, forcing $t=2$. Hence, $N=\PGL(2,q_0)$ with $q=2^f=q_0^2$.
It follows from $2^{f-x}+1\mid3$ that $f=x+1$ and $k=3$, implying $|B_1\cup B_2\cup B_3|=5$.
Note that $X_{\alpha,\beta}$ has orbit sizes 1, 1 and $q-1$,
and that $B_1\cup B_2\cup B_3$ is a union of $G_{\alpha,\beta}$-orbits.
We have that $|B_1\cup B_2\cup B_3|=q+1$, and so $q=4$, $(v,b,r,k)=(5,10,6,3)$. 
Since $\PSL(2,4)\unlhd G \leq \mathrm{P \Gamma L}(2,4)$ acts $2$-transitively on $\mathcal{D}$, it follows that $\mathcal{D}$ is the complete $2$-$(5,3,3)$ design.

Suppose that $N=C_p^f \rtimes C_{q-1/\gcd(2,q-1)}$. Then $|X:N|=q+1=v$ divides $b$,
and hence from $vr=bk$ we have $k\mid r$.
Therefore, $p^{f-x}+1 \mid 3p^x$, and so $p^{f-x}+1=3$, implying
$p=2$, $f=x+1$ and $k=3$. By the same argument as above, 
we have that $\mathcal D$ is a complete $2$-$(5,3,3)$ design.

Suppose finally that $N=\mathrm{D}_{2(q-\epsilon)/\gcd(2,q-1)}$ with $\epsilon=\pm1$. 
Then 
\[
q(q+\epsilon)(p^{f-x}+1)\mid 6p^x(q+1),
\]
and so $q\mid 6p^x$, implying $f=x+1$ and $p=2$ or 3.
If $p=3$, then $2(3^f+\epsilon)\mid 3^f+1$, forcing $\epsilon=-1$ and $f=1$, a contradiction.  
Therefore, $p=2$, $f=x+1$ and $k=3$. Again, we have that $\mathcal D$ is a complete $2$-$(5,3,3)$ design.
\end{proof}

\textbf{Acknowledgements.}

The first author's work is supported by the National Natural Science Foundation of China (No:12101120).

All of the authors thank the Italian National Group for Algebraic and Geometric Structures and their Applications (GNSAGA–INdAM) for its support to their research.

\end{document}